\magnification=\magstep 1
\baselineskip=15pt     
\parskip=3pt plus1pt minus.5pt
\overfullrule=0pt
\font\hd=cmbx10 scaled\magstep1

\def\Pone{{\bf P}^1}
\def\P{{\bf P}}

\def\C{{\cal C}}

\def\O{{\cal O}}

\def\H{{\cal H}}
\def\G{{\cal G}}

\def\F{{\cal F}}

\def\M{{\cal M}}

\def\E{{\cal E}}

\def\EY{\E(1)|_Y}

\def\PEY{\P(\EY)}

\def\OP{\O_{\PEY}(1)}
\def\OP1{\O_{\Pone}}

\def\Pic{{\rm Pic}}

\def\re{\mathop{\rm Re}}
\def\im{\mathop{\rm Im}}

\def\boldz{{\bf Z}}

\def\dual#1{{#1}^{\scriptscriptstyle \vee}}

\def\mapright#1{\smash{
  \mathop{\longrightarrow}\limits^{#1}}}

\def\mapdown#1{\Big\downarrow
   \rlap{$\vcenter{\hbox{$\scriptstyle#1$}}$}}

\input epsf.tex
\input amssym.tex
\centerline{\hd Mirror Symmetry via 3-tori for a class of 
Calabi-Yau Threefolds}
\medskip
\centerline{\it Mark Gross\footnote{*}{Supported in part by NSF grant 
DMS-9400873 and Trinity College, Cambridge}\hfil P.M.H. Wilson}
\medskip
\centerline{November 8th, 1996}
\medskip
{\settabs 3 \columns
\+Department of Mathematics&&Department of Pure Mathematics\cr
\+Cornell University&&University of Cambridge\cr
\+Ithaca, NY 14853&&Cambridge CB2 1SB, U.K.\cr
\+mgross@math.cornell.edu&&pmhw@pmms.cam.ac.uk\cr}

\bigskip
\bigskip
{\hd \S 0. Introduction}

The authors of [16] have proposed a conjectural construction of mirror
symmetry for Calabi-Yau threefolds. They argue from the physics that in a 
neighbourhood of the large complex structure limit (see [11] for the 
definition of large complex structure limits), any Calabi-Yau threefold
$X$ with
a mirror $Y$ should admit a family of supersymmetric toroidal 3-cycles.
In mathematical terminology, this says that there should be a fibration
on $X$ whose general fibre is a special Lagrangian 3-torus $T^3$.

We recall from [8] the notion of a special Lagrangian submanifold.

\proclaim Definition.  Let $X$ be a K\"ahler manifold of dimension $n$,
with complex structure $I$ and K\"ahler metric $g$, whose holonomy is
contained in $SU(n)$.  Recall that this latter condition is equivalent
to the existence of a covariant constant holomorphic $n$-form $\Omega$
on $X$.  We say that $\Omega$ is normalized if
$$ (-1)^{n(n-1)/2} (i/2)^n \Omega \wedge \bar\Omega = \omega ^n /n! $$
the volume form on $X$, where $\omega$ denotes the K\"ahler form -- 
thus the normalized holomorphic $n$-form on $X$ is unique up to a phase
factor $e^{i\theta}$.  We say that a submanifold $M$ of $X$ is special 
Lagrangian (we shall take $M$ to be embedded, although more generally 
one only needs it to be immersed) if a normalized holomorphic $n$-form 
$\Omega$ can be chosen on $X$ whose real part $\re \Omega$ restricts to
the volume form on $M$ -- the normalization has been chosen here 
so that $\re \Omega$ is a calibration in the sense of [8].  This 
definition is equivalent to $M$ being an oriented submanifold of real
dimension $n$ such that the restriction of $\omega$ to $M$ is zero,
i.e. $M$ is Lagrangian, together with the existence of a holomorphic 
$n$-form $\Omega$ whose imaginary part $\im \Omega$ also restricts
to zero on $M$ [8].

Given a special Lagrangian submanifold $M$ of $X$, it is shown in [10]
that the deformations of $M$ as a special Lagrangian submanifold
are unobstructed and that the tangent space to the local
deformation space may be identified with the harmonic 1-forms on $M$.
The dimension of the local deformation space is therefore
$\dim_{\bf R} H^1(M,{\bf R})$. For $X$ a Calabi-Yau $n$-fold and $T$ a special
Lagrangian $n$-torus submanifold, 
$T$ therefore moves in an $n$-dimensional family.

\proclaim Definition. A Calabi-Yau $n$-fold $X$ is said to have a
special Lagrangian $n$-torus fibration if there exists a map of
topological manifolds $f:X\rightarrow B$ whose general fibre is a special
Lagrangian $n$-torus. 

We have taken $B$ as a topological manifold, but
from the results of [10], it will be locally a differentiable 
manifold with a natural Riemannian metric [10], (3.10), 
except perhaps at points corresponding to singular fibres.
In \S 3, we shall generalise the notion of special Lagrangian $n$-torus
fibration to the case when the metric on $X$ is allowed to degenerate
in a suitably nice way.

Let us now take $X$ to be a Calabi-Yau threefold. If $X$ contains
a special Lagrangian $T^3$, then it moves locally in a three dimensional
family. There are hard questions which need to be addressed
concerning whether the family obtained foliates the manifold and whether
it can be suitably compactified so as to obtain a special Lagrangian
3-torus fibration on $X$--- for some discussion of these
problems and further discussion of the motivation from physics,
we refer the reader to [12], where an account of [16] for
mathematicians will be found. In this paper, we study a class of
Calabi-Yau threefolds for which these issues are not a problem (provided
we work with a mildly degenerate metric), and $X$ does have a special
Lagrangian torus fibration.

Under this assumption, [16] provides a recipe for constructing the 
mirror $\check X$ to $X$. Topologically $\check X$ is a compactification of
the moduli space parametrising special Lagrangian 3-tori (with fixed
cohomology class) on $X$, together with a flat $U(1)$-connection
on $T$. For a given torus $T$, the flat $U(1)$-connections
are parametrized by the dual torus, so the recipe may be rephrased
as follows. We take the smooth part of the special Lagrangian $n$-torus
fibration, which we denote by $p:X_0\rightarrow B_0$. The fibres of $p$
are therefore all special Lagrangian tori. We then dualize the fibres;
explicitly we form the family of tori $\dual{X_0}$ over $B_0$ given by
$R^1p_*{\bf R}/R^1p_*\boldz$. The mirror
$\check X$ (which has been assumed to exist) is recovered as a 
compactification of this dual special Lagrangian 3-torus fibration.
We observe that if $X$ and $\check X$ are to be mirrors
in the usual sense (in particular the mirror of $\check X$ is $X$),
then the special Lagrangian $3$-torus fibration on $X$ should have
a special Lagrangian section. There is also an argument
from the physics that this is the case, and it is true for the
examples we study below provided we restrict attention to suitable subsets
of moduli.
Note here that the Euler characteristic $e(X)$ is concentrated purely
in the singular fibres of the 3-torus fibration, and so we expect to see
locally how the Euler characteristic changes sign when we pass
to the mirror. The authors of [16] also suggest a method for putting
a complex structure on the dual fibration over $B_0$--- the
difficulties here are also discussed in [12].

We shall refer to the recipe from [16] as the SYZ construction.
The current paper explores the SYZ construction in two cases. The mirror
map for K3 surfaces has been studied by several authors, and so
the recipe from [16] should reproduce known results. This is checked in
\S 1 and the beginning of \S 4. A different approach to this
appears in [12]. In the case of Calabi-Yau threefolds,
the conjectured construction has not previously been worked out for
any examples. Most of this paper will therefore be devoted to
doing this
for a class of Calabi-Yau threefolds
studied independently by Ciprian Borcea and Claire Voisin.
These are obtained from a K3
surface $S$ with holomorphic involution $\iota$,
by resolving singularities of the 
threefold $S\times A/(\iota,j)$, where $A$ is an elliptic curve and $j$ 
is the involution given by negation. It is known when the mirror family exists
for such a Calabi-Yau threefold, and when it exists there is a very natural
description of it [3, 17], the construction largely depending on
ideas of Nikulin [13]. In the case of these Borcea-Voisin threefolds,
we do have natural choices for the special Lagrangian torus
fibration, and we check (3.1) that the mirror does, at least topologically,
satisfy the properties required by the SYZ construction. The question
of complex and K\"ahler structures (i.e. the mirror map)
is considered further in \S 4. But even just at the topological
level, there is a very beautiful description of how passing to the mirror
affects the singular fibres of the special Lagrangian torus
fibration and how this causes the Euler characteristic to change sign.

In the final section, we discuss the mirror map (from the point of view
of special Lagrangian torus fibrations) for both the K3 case and the
Borcea-Voisin examples. More generally, these examples lead us to
a conjectural interpretation
of the mirror map for Calabi-Yau threefolds whose mirrors can be constructed
by the SYZ method, in particular having special Lagrangian 3-torus
fibrations. This conjectural construction involves the Leray spectral
sequence associated to the $3$-torus fibration, and is consistent
with the results obtained for the Borcea-Voisin examples.

The authors with to thank Nigel Hitchin and David Morrison for
very helpful discussions during the preparation of this manuscript.

\bigskip
\vfill
\eject
{\hd \S 1. Mirror Symmetry for K3 Surfaces}

We recall here the construction of mirror symmetry for K3 surfaces
given in 
[6]. This construction was discovered independently by Pinkham [15]
and Dolgachev and Nikulin [7] (see [6] for precise references
and more history). This construction also is a special case of a more
general construction due to Aspinwall and Morrison [1] which is
more directly inspired by the physics involved. It is possible
to interpret this more general form of mirror symmetry in terms of the
SYZ construction. This is done in [12].

Let $L$ be the K3 lattice, i.e. an even unimodular lattice of signature
$(3,19)$, isomorphic to $U(1)\oplus U(1)\oplus U(1)\oplus E_8\oplus
E_8$, where $U(m)$ is the rank 2 lattice with intersection matrix
$\pmatrix{0&m\cr m&0\cr}$. If $S$ is a K3 surface then $H^2(S,\boldz)
\cong L$.

Now suppose that $S$ is a K3 surface, and fix a marking $\phi:
H^2(S,\boldz)\mapright{\cong} L$. We set
$$\eqalign{M&=\Pic(S)\subseteq L\cr
T&=M^{\perp}\subseteq L.\cr}$$

\proclaim Definition. An isotropic vector $E\in T$ is called {\it $m$-admissible
} if there exists an isotropic vector $E'\in T$ such that $E.E'=m$,
and there does not exist a vector $\alpha\in T$ with either
$0<\alpha.E<m$ or $0<\alpha.E'<m$.

By [6], Lemma (5.4), this is equivalent to the existence of
an isotropic vector $E'\in T$ with the sublattice $P$ of $T$
generated by $E$ and $E'$ isomorphic to $U(m)$, and $T$ has an orthogonal
decomposition $P\oplus P^{\perp}$. 

If $E\in T$ is an $m$-admissible vector, then we define
$$\check M=(\boldz E)_T^{\perp}/\boldz E\subseteq T/\boldz E.$$
There is a natural primitive embedding $i:\check M\rightarrow T$
given by 
$$i(\alpha)=\alpha-{\alpha.E'\over m} E.$$
In fact $T=P\oplus\check M$.

Now to a given primitive sublattice $M\subseteq L$ of signature
$(1,t)$, $T=M^{\perp}$, one can associate to it the period domain
of {\it marked $M$-polarized K3 surfaces}
$$D_M=\{{\bf C}\Omega\in \P(T\otimes_{\boldz}{\bf C})\ |\ 
\Omega.\Omega=0, \Omega.\bar\Omega>0\}.$$
By the Torelli theorem for K3 surfaces, $D_M$ is the moduli space of
marked K3 surfaces $(S,\phi)$ with $\phi:H^2(S,\boldz)\mapright{\cong}
L$ a marking such that $\phi^{-1}(M)\subseteq\Pic S$. Note that $D_M$
has two connected components interchanged by complex conjugation.
The moduli space of $M$-polarized K3 surfaces will be a suitable
quotient of $D_M$ via a group action, but we will not be concerned about
this action here and will work on the level of period domains. See [6]
for more details.

We consider also the tube domain 
$$T_M=\{B+i\omega\in M\otimes_{\boldz}{\bf C}\ |\ \omega.\omega>0\}.$$
The class $B$ above is referred to by physicists as the $B$-field,
but usually represents a class in $H^2(S,{\bf R}/\boldz)$. The tube
domain $T_M$ represents a cover of the complex
K\"ahler moduli space, analogous to the period domain being a cover
of the complex moduli space. We shall always construct the mirror
map (which gives an isomorphism between K\"ahler moduli and complex
moduli of the mirror) at the level of these covers.

We have a slightly more explicit version of [6], Theorem (4.2).

\proclaim Proposition 1.1. Given a choice of splitting $T=P\oplus \check M$,
$P\cong U(m)$, $E,E'\in P$ primitive isotropic vectors with $E.E'=m$,
there is an isomorphism $$\phi:T_{\check M}\rightarrow D_M$$
given by 
$$\eqalign{\phi(\check B+i\check\omega)= &\check B+{E'\over m}
+\left({\check\omega.\check\omega-\check B.\check B\over 2}\right) E\cr
&+i(\check \omega-(\check\omega.\check B)E)\in \P(T\otimes_{\boldz}{\bf C}).
\cr}$$

Proof: If $\Omega=\phi(\check B+i\check \omega)$, it is easy to check
that $\Omega.\Omega=0$ and $\Omega.\bar\Omega=2\check\omega.\check\omega>0$.
Now by [6], Lemma (4.1), if ${\bf C}\Omega\in D_M$, then $\Omega.E\not=0$.
Thus we can normalize $\Omega$ by multiplying by a complex constant
to ensure that $\Omega.E=1$. The map $\phi^{-1}$ can then be defined
by taking $\phi^{-1}({\bf C}\Omega)$ to be $\check B+i\check \omega$,
where $\check B$ and $\check \omega$ are the orthogonal projections
onto $\check M$ of $\re\Omega$ and $\im\Omega$ respectively, using
the orthogonal decomposition $T=P\oplus\check M$. $\bullet$

The case that $m=1$ is especially important; in this case, given 
$P\cong U(1)\subseteq T$, we also have $\check M^{\perp}\cong P\oplus M$,
and so repeating the process again we recover $M$. 
In this case, we get isomorphisms
$$\check\phi:T_M\rightarrow D_{\check M}$$
and 
$$\phi:T_{\check M}\rightarrow D_M.$$
We say in this situation that the family of $M$-polarized K3 surfaces
and $\check M$-polarized K3 surfaces are {\it mirror families}.
$D_{\check M}$ and $D_M$ parametrize the mirror families of marked K3 surfaces,
and $\check \phi$ and $\phi$ are the mirror maps exchanging K\"ahler
and complex moduli. So in particular, the family of $M$-polarized K3 surfaces
has a mirror family if and only if $T$ contains a sublattice isomorphic
to $U(1)$.

We now describe how this picture coincides with the proposed construction
of mirror symmetry in [16]. We first need to understand which 
submanifolds of a K3 surface are special Lagrangian.

Fix a K3 surface $S$ with complex structure $I$ and K\"ahler-Einstein
metric $g$. $(S,g)$ is a hyperk\"ahler manifold with complex structures
$I$, $J$ and $K$ generating an $S^2$ of possible complex structures for
which $g$ is a K\"ahler metric. The period of $S$ in complex structure
$I$ is given by the complex 2-form $\Omega$ with
$$\Omega(X,Y)=g(JX,Y)+ig(KX,Y),$$
and the K\"ahler form on $S$ is given by
$$\omega(X,Y)=g(IX,Y).$$
Thus we have, in the various complex structures $I$, $J$ and $K$, the
following data:
\bigskip
{\settabs 3\columns
\+Complex Structure&Holomorphic 2-form&K\"ahler form\cr
\+$I$&$\re\Omega+i\im\Omega$&$\omega$\cr
\+$J$&$\omega+i\re\Omega$&$\im\Omega$\cr
\+$K$&$\im\Omega+i\omega$&$\re\Omega$\cr

}
Here $$(\re\Omega)^2=(\im\Omega)^2=\omega^2>0.$$
We observe that $\Omega$ is then 
normalized, in the sense that was defined in the Introduction.
Note that fixing the metric $g$ and complex structure $I$ still
allows an $S^1$ of choices for $J$ and $K$, which corresponds to
multiplying $\Omega$ by a phase $e^{i\theta}$.

The following observation is due to Harvey and Lawson ([8], p. 154).

\proclaim Proposition 1.2. A two-dimensional submanifold $Y\subseteq S$
is special Lagrangian with respect to $(g,I)$ if it is holomorphic
with respect to a complex structure $K$, where $K$ is one
of the possible complex structures mentioned above.

Proof: $Y\subseteq S$ is special Lagrangian with respect to $(g,I)$
if $(\re\Omega)|_Y$ is the volume form on $Y$ induced by the 
metric $g$ for some suitable choice of phase of $\Omega$. But this
is equivalent to $Y$ being holomorphic in complex structure $K$
by Wirtinger's theorem (see [8], p. 58, Example I). $\bullet$

We can now make the connection between the existence
of a mirror family to the family $D_M$ of $M$-polarized K3 surfaces
and the existence of a special Lagrangian torus fibration on $S$
for general choice of complex structure on $S$ in $D_M$. 

\proclaim Proposition 1.3. Let $M\subseteq L$ be a primitive
lattice of signature $(1,t)$ for some $t$. 
Suppose there is an $m$-admissible vector $E\in T$.
Then for general choice of K3 surface $S$ in $D_M$ and
for general choice of K\"ahler-Einstein metric $g$ on $S$ compatible
with its complex structure, $S$ has a fibration $\pi:S\rightarrow S^2$ (the
two-sphere)
whose general fibres are special Lagrangian 2-tori, and such that $\pi$
has a special Lagrangian numerical $m$-section.
Here, a {\it numerical $m$-section}
is a surface in $S$ whose topological intersection number with a fibre
of $\pi$ is $m$.

Proof: If $E\in T$ is $m$-admissible, we obtain
a decomposition $T=P\oplus\check M$ with $P\cong U(m)$. If we choose
$\omega\in M\otimes_{\boldz}{\bf R}$ to be a general K\"ahler class,
then we can assume that
$\omega^{\perp}\cap M=0$. We also choose the complex structure
on $S$, but initially not a general one: choose a complex structure on $S$
so that its holomorphic 2-form is $\Omega=\phi(\check B+i\check\omega)$
for a general choice of $\check B+i\check \omega\in T_{\check M}$ such
that $\check B.\check \omega=0$. From the formula
of Proposition 1.1, it follows that 
$$E.(\im\Omega+i\omega)=E'.(\im\Omega+i\omega)=0.$$
Thus if we denote by $S_K$ the K3 surface with the corresponding complex
structure $K$, then $U(m)\subseteq \Pic S_K$
and $\Pic S_K\subseteq T$. Thus
$\pm E$ is an effective class
on $S_K$. Without loss of generality, we can assume $E$ is effective, 
and after reflecting $E$ around some $-2$ curves in $\Pic S_K$ then
by [14], \S 3 Cor. 3 and \S 6 Theorem 1, we can
assume furthermore that $E$ is the class of a fibre of an elliptic
fibration $f:S_K\rightarrow\Pone\cong S^2$. 
This yields a special
Lagrangian fibration $\pi : S \to S^2$ with respect to $(g,I)$, 
where $g$ is the Ricci-flat metric corresponding to $\omega$.
Furthermore, since $E'$ is algebraic on $S_K$, a component of
$\pm E'$ will be a holomorphic $m$-section of $f$, and
hence will yield a special Lagrangian $m$-section of $\pi$.

This proves the implication holds for general choice of
K\"ahler class and somewhat special choice of complex structure.
To extend the result to general choice of complex
structure, let $\sigma$ be the cohomology class of the special
Lagrangian $m$-section constructed above. For general choice
of $\Omega\in D_M$ on $S$, and some suitable choice of phase
for $\Omega$, $\sigma$ will be algebraic on $S_K$. Now $\sigma$ being
represented by an irreducible curve is an open condition on the 
moduli of K3 surfaces in which $\sigma$ is algebraic, and also
we have seen above that for certain values
of the complex structure, $\sigma$ is represented
by an irreducible curve. It then follows that for general
choice of complex structure on $S$, $\sigma$ is represented by an
irreducible curve on $S_K$. This then gives a special Lagrangian
submanifold $\sigma$ on $S$, with $\sigma.E=m$. 

On the other hand, if we suitably change the phase of $\Omega$,
then $E$ is algebraic on a surface $S_{K'}$ (not necessarily the
same surface as $S_K$ above), and as above, this yields a 
special Lagrangian torus fibration $\pi:S\rightarrow S^2$. The special
Lagrangian submanifold $\sigma$ constructed above is then a numerical
$m$-section of $\pi$.
$\bullet$

{\it Remark.} In the above proposition, we need to insist on ``numerical''
since as the curves $E$ and $\sigma$ may not be algebraic with respect to
the same complex structure, there is no guarantee that all the intersection
multiplicities are positive.
However, as we saw in the proof of the above
Proposition, there are certain complex structures where both $E$ and
$\sigma$ are algebraic on $S_K$ for the same complex structure $K$. 
In this case, the numerical $m$-section is a genuine $m$-section;
in fact, it is clear that the special Lagrangian numerical
$m$-section we found above continues to be a genuine $m$-section
on some open (but not necessarily Zariski open or dense) subset of $D_M$.

\proclaim Corollary 1.4. The family of $M$-polarized K3 surfaces has
a mirror family if and only if for a general K3 surface $S$ with $\Pic S\cong
M$, and for any general choice of K\"ahler metric on $S$,
$S$ has a special Lagrangian 2-torus fibration with a special Lagrangian
numerical section.

Thus we see that, as claimed in [16], the existence of a mirror
of a K3 surface is equivalent to the (generic) existence of a special
Lagrangian fibration with a special Lagrangian (numerical) section.

To produce the mirror topologically, we must choose $S$ in the family for which
the 2-torus fibration has a special Lagrangian section and then dualize the torus
fibration. As is well-known, this does not change the topology
of the surface, but we wish to make this explicit here.

Let $B_0\subseteq \Pone({\bf C})$ be the locus over which the fibration
$f$ is smooth, and let $S_0=f^{-1}(B_0)$. Let $p:S_0\rightarrow
B_0$ be the restriction of $f$. Then $\dual{S_0}\rightarrow B_0$
is $R^1p_*{\bf R}/R^1p_*\boldz\rightarrow B_0$, and $S_0\rightarrow
B_0$ can be identified (only topologically) with
$\dual{(R^1p_*{\bf R})}/\dual{(R^1p_*\boldz)}\rightarrow B_0$. 
Poincar\'e duality gives a perfect pairing
$$R^1p_*{\bf R}\times R^1p_*{\bf R}\rightarrow R^2p_*{\bf R}\cong {\bf R}$$
and thus gives a natural isomorphism $R^1p_*{\bf R}\cong\dual{(R^1p_*{\bf R})}$.
Using this isomorphism, we identify $S_0\rightarrow B_0$ with $\dual{S_0}
\rightarrow B_0$, and so the original fibration $S_K\rightarrow \Pone({\bf C})$
may be identified topologically as a compactification of the dual fibration
$\dual{S_0} \rightarrow B_0$.

We will use this identification in a crucial way in \S 2. However,
as noted in [16], we have only found that topologically the mirror
of a K3 surface is a K3 surface; we also need to explain the mirror map.
We take this up in \S 4.

\bigskip

{\hd \S 2.  K3 Surfaces with Involution}

Now let $S$ be a K3 surface equipped with an involution $\iota:S\rightarrow
S$, acting by $-1$ on the holomorphic 2-form.
Let $C_1,\ldots,C_N$ be the fixed curves of $\iota$. By [17], (1.1), these
are smooth curves and either
$$\vbox{\item{(1)} $N=0$ (in which case $\iota$ is the Enriques involution).
\item{(2)} $N=2$ and $C_1$ and $C_2$ are elliptic curves.
\item{(3)} $C_2,\ldots,C_N$ are rational, and $C_1$ has genus
$N'\ge 0$.}$$

$\iota$ induces a map on cohomology which we will write as
$H(\iota):H^2(S,\boldz)\rightarrow H^2(S,\boldz)$. Using the same notation
as in \S 1, we have $M=H^2(S,\boldz)^+$, the group of invariants of 
$H(\iota)$, and $T=H^2(S,\boldz)^-$, the group of anti-invariants
of $H(\iota)$. Let $\Delta\subseteq D_M$ be defined by
$$\Delta=\{{\bf C}\Omega\in D_M\ |\ 
\hbox{there exists an $\alpha\in T$, $\alpha\not=0$ such that
$\alpha.\Omega=0$.}\}$$
Then $D_M\backslash \Delta$ is the period domain for marked K3 surfaces with
involution acting on cohomology via $H(\iota)$. (See [17] (2.1.1).)

Suppose furthermore that there exists a $P\subseteq T$, $P\cong U(1)$. Then
$T$ decomposes as $T=P\oplus\check M$.
Following [17] we define $r_P:H^2(S,\boldz)\rightarrow H^2(S,\boldz)$
by $r_P|_P=Id_P$, $r_P|_{P^{\perp}}=-Id_{P^{\perp}}$. Then define
$H(\check \iota):H^2(S,\boldz)\rightarrow H^2(S,\boldz)$ by $H(\check 
\iota)=r_P\circ
H(\iota)$. This will induce an involution $\check \iota$ on each member
$\check S$ of the mirror family of $\check M$-polarized K3 surfaces,
by Torelli. It was observed by Borcea [3] and
Voisin [17] that the work of Nikulin
[13] implies that
passing from $\iota$ to $\check\iota$ interchanges the numbers
$N$ and $N'$ in Case (3) above.
After fixing a basis $E,E'$ for $P$, the mirror maps
$$\check \phi:\check\phi^{-1}(D_{\check M}\backslash \check\Delta)
\rightarrow D_{\check M}\backslash \check \Delta$$
and
$$\phi:\phi^{-1}(D_{M}\backslash \Delta)
\rightarrow D_{M}\backslash  \Delta$$
give maps between the moduli of K3s with involution of a given topological
type on the one hand and an open subset of the tube domain of the mirror.

We will now show how the process of dualising special Lagrangian 2-torus
fibrations will produce the mirror involution. To do this, we
make precise the natural notion of a dual involution.

Let $\pi:X\rightarrow B$ be a smooth $n$-torus fibration with
a distinguished section, so as to make $\pi$ into a fibration of topological
groups. (This is the same as choosing an isomorphism of $\pi$
with
the $n$-torus fibration $\dual{(R^1\pi_*{\bf R})}/\dual{(R^1\pi_*\boldz)}$.)
A {\it homomorphism} of smooth $n$-torus fibrations $(\pi_1:X_1\rightarrow
B_1)\mapright{(\alpha,\beta)} (\pi_2:X_2\rightarrow B_2)$ will then be 
a commutative diagram
$$\matrix{X_1&\mapright{\alpha}&X_2\cr
\mapdown{\pi_1}&&\mapdown{\pi_2}\cr
B_1&\mapright{\beta}&B_2\cr}$$
such that the map $X_1\rightarrow X_2\times_{B_2} B_1$ is 
induced by an ${\bf R}$-linear
map
$\alpha:
\dual{(R^1\pi_{1*}{\bf R})}
\rightarrow
\beta^*\dual{(R^1\pi_{2*}{\bf R})}$.
If $(\alpha,\beta)$ 
is an isomorphism, then we can naturally define the transpose
of $(\alpha,\beta)$, which is a map
$$(\alpha,\beta)^t:
(\dual{\pi_2}:\dual{X_2}\rightarrow B_2)
\rightarrow
(\dual{\pi_1}:\dual{X_1}\rightarrow B_1)
$$
given by the commutative diagram
$$\matrix{\dual{X_2}&\mapright{\alpha^t}&\dual{X_1}\cr
\mapdown{\dual{\pi_2}}&&\mapdown{\dual{\pi_1}}\cr
B_2&\mapright{\beta^{-1}}&B_1\cr}$$
where $\alpha^t$ is induced by the induced map
$$\alpha^*:R^1\pi_{2*}{\bf R}\rightarrow (\beta^{-1})^*R^1\pi_{1*}{\bf R}.$$

We now show that for certain values of complex structure on $S$ and
$\check S$, the transpose
involution is precisely the mirror involution. 
We first consider a subset $D_M'\subseteq D_M$ given
by
$$D_M'=\{{\bf C}\Omega\in D_M\ |\ \im\Omega\in \check M\otimes_{\boldz}
{\bf R}\}.$$
Note that this is a {\it real} codimension one condition. We define
$D'_{\check M}$ similarly. Put also
$$T_M'=\{B+i\omega\in T_M\ |\ B.\omega=0\}.$$
$T_M'$ is a real codimension one subset of $T_M$. We define $T'_{\check M}$
similarly.

We now fix a marked K3 surface with involution $(S,\iota)$ whose period
${\bf C}\Omega\in D_M'\backslash \Delta_M$, 
and fix a K\"ahler form $\omega$ on $S$ normalized
so that 
$(\re\Omega)^2=(\im\Omega)^2=\omega^2$. This determines a K\"ahler-Einstein 
metric $g$ on $S$ , and $\Omega$ is then normalized with respect
to the pair $(g,I)$.
Let $S_K$ be the K3 surface with complex
structure $K$. $S_K$ has period $\Omega_K=\im \Omega+i\omega$,
so it is clear that, for general choice of period ${\bf C}\Omega$ in 
$D_M'\backslash \Delta_M$ and general choice of $\omega$, $\Pic S_K
=P\cong U(1)$. Thus $S_K$ has a Jacobian elliptic fibration $f:S_K
\rightarrow \Pone({\bf C})$, and we can choose a basis $E$, $E'$ for 
$P$ such that $E$ is the class of a fibre of $f$ and $E'-E$ is the
class of the (unique) section of $f$. Thus, roughly speaking,
$D_M'$ is the period space for marked $M$-polarized K3 surfaces $S$ such that
the special Lagrangian 2-torus fibration on $S$ and the special Lagrangian
section of this fibration are holomorphic with respect to the same
complex structure $K$. In particular, we know the fibration has
a special Lagrangian section, not just a numerical section.

Having fixed $E$ and $E'$, we have now fixed the mirror maps
$\check \phi:T_M\rightarrow D_{\check M}$, $\phi:T_{\check M}\rightarrow
D_M$. These have the property that $\check\phi(T_M')=D_{\check M}'$
and $\phi(T'_{\check M})=D_M'$, as can be easily checked.
Now choose $\check S$ which is mirror to $S$
to have period ${\bf C}\check\Omega=\check 
\phi(B+i\omega)$ for any choice of $B$ for which $B+i\omega\in
T_M'\backslash \check\phi^{-1}(\Delta_{\check M})$ and let
$\check B+i\check \omega=\phi^{-1}({\bf C}\Omega)$. Thus we can write
$$\Omega=\phi(\check B+i\check \omega)=\check B+E'+\left({\check 
\omega.\check \omega-\check B.\check B\over 2}\right )
E+i\check \omega$$
and
$$\check\Omega=\check\phi(B+i\omega)=
B+E'+\left({\omega.\omega-B.B\over 2}\right )
E+i\omega.$$
Note that $\Omega \wedge \bar\Omega = 2 \check\omega ^2 $
and the metric $g$ has been chosen so that $\Omega \wedge \bar\Omega =
2 \omega ^2 $, and hence $\check \omega ^2 = \omega ^2.$

Let $\check g$ be the K\"ahler-Einstein metric on $\check S$ compatible
with $\check\omega$ and the complex structure $\check I$ on $\check S$,
and let $\check J$, $\check K$ be the complex structures on $\check S$
for which $\check\Omega(X,Y)=\check g(\check J X,Y)+i\check g(\check K X,Y)$.

To summarize, we have the following table of normalized holomorphic 2-forms and
K\"ahler forms:

{\settabs=3 \columns
\+Complex structure&Holomorphic 2-form&K\"ahler form\cr
\+$I$&$\Omega$&$\omega$\cr
\+$K$&$\check\omega+i\omega$&$\check B+E'+\left(
{\check\omega.\check\omega-\check B.\check B\over 2}\right )E$\cr
\+$\check I$&$\check \Omega$&$\check \omega$\cr
\+$\check K$&$\omega+i\check \omega$&$ B+E'+\left(
{\omega.\omega- B. B\over 2}\right )E$\cr

}

The mirror pair of 
K3 surfaces with involution $(S,\iota)$ and $(\check S,\check\iota)$
will now have a direct topological relationship obtained by dualising
the 2-torus fibration. More precisely

\proclaim Theorem 2.1. Let $B_0\subseteq \Pone({\bf C})$ be the locus over
which $f:S_K\rightarrow \Pone({\bf C})$ is smooth, and let $S_0=f^{-1}(B_0)$,
thinking of $S_0$ only as a topological manifold with 2-torus
fibration $S_0\rightarrow B_0$. Then $\iota$ and $\check\iota$ induce
homomorphisms of 2-torus fibrations
$$\hbox{
$\matrix{S_0&\mapright{\iota|_{S_0}}&S_0\cr
\mapdown{}&&\mapdown{}\cr
B_0&\mapright{\iota'}&B_0\cr}\quad\quad$
and 
$\quad\quad\matrix{S_0&\mapright{\check\iota|_{S_0}}&S_0\cr
\mapdown{}&&\mapdown{}\cr
B_0&\mapright{\check\iota'}&B_0\cr}$}$$
Furthermore, $\iota'=\check\iota'$ and 
if $\Psi:S_0\rightarrow \dual{S_0}$ is the isomorphism given by
Poincar\'e duality as in \S 1, then $\Psi^{-1}\circ
(\iota|_{S_0})^t\circ \Psi=\check\iota|_{S_0}$.

{\it Remark.} This indicates a very special relationship between the pairs
$(S,\iota)$ and $(\check S,\check\iota)$ which will only exist for very
particular choices of complex structure on $S$ and $\check S$. Indeed, as we 
change the complex structure on $\check S$, we would expect the map $\check
\iota$ to change when thought of as a map of topological manifolds, even
though the homotopy class of $\check\iota$ does not change. We remark
that choosing complex structure in $D_M'$ and K\"ahler structure in
$T'_M$ enables us to simplify considerably the ensuing analysis.
In particular,
the careful choice of complex and K\"ahler structures yields the following
additional structure for the maps $\iota$ and $\check\iota$.

\proclaim Lemma 2.2. As maps on $S_K$, $\iota$ and $\check\iota$ are
anti-holomorphic. 

Proof: We first show that $\iota:S_K\rightarrow S_K$ is anti-holomorphic.
The metric $g$ is also the K\"ahler-Einstein
metric on $S_K$. Thus $\iota$ is an isometry
with respect to this metric on $S_K$, and takes harmonic forms to
harmonic forms. If $\Omega_K$ is the holmorphic 2-form on $S_K$, whose
de Rham cohomology class is $\check \omega+i\omega$, then $\Omega_K$ is
harmonic, and thus so is $\iota^*\Omega_K=-\check\omega+i\omega$ 
in cohomology. Thus
$\iota^*\Omega_K=-\bar\Omega_K$ as a two-form, not just as a cohomology
class. From this we see that 
if locally $\Omega_K$ is of the form $f dz_1\wedge dz_2$,
$f$ holomorphic, then $\iota^*\Omega_K$ is of the form $-\bar f d\bar z_1\wedge
d\bar z_2$. Thus the Jacobian of $\iota$ with respect to the variables
$z_i,\bar z_i$ must be of the form $$\pmatrix{0&0&*&*\cr 0&0&*&*\cr 
*&*&0&0&\cr
*&*&0&0&\cr}$$
and so $\iota$ is anti-holomorphic.

The same proof shows that $\check\iota:\check S_{\check K}\rightarrow\check
S_{\check K}$ is anti-holomorphic. Since $\check K$ and $K$ have
complex conjugate periods, $\check\iota:S_K\rightarrow S_K$
is also anti-holomorphic.
$\bullet$

{\it Proof of Theorem 2.1.} Since $\iota:S_K\rightarrow S_K$
is anti-holomorphic and takes the class $E$ to $-E$
and the class of the unique section $\sigma$ to $-\sigma$, $\iota$
preserves the fibration and the section of $f:S_K\rightarrow \Pone$; i.e.
we have a commutative diagram
$$\matrix{S_K&\mapright{\iota}&S_K\cr
\mapdown{f}&&\mapdown{f}\cr
\Pone({\bf C})&\mapright{\iota'}&\Pone({\bf C})\cr}$$
with $\iota'$ the anti-holomorphic involution on $\Pone({\bf C})$
induced by identifying $\Pone({\bf C})$ with the section $\sigma$. In 
particular, restricting this diagram to $S_0\rightarrow B_0$,
we see that $\iota|_{S_0}$ yields a homomorphism of 2-torus
fibrations. Here $\iota|_{S_0}$ being a homomorphism follows from
the anti-holomorphicity of $\iota$ in exactly the same
way it would follow if $\iota$ was holomorphic. A similar argument works
for $\check \iota$, and this proves the first part of Theorem 2.1.

Let $-1:S_K\rightarrow S_K$ be the holomorphic involution given by fibrewise
negation of $f:S_K\rightarrow \Pone({\bf C})$. We now finish the proof by 
verifying the following two claims:

{\it Claim 1:} $\check\iota=(-1)\circ \iota$.

{\it Claim 2:} Under the identification of $S_0\rightarrow B_0$
with $\dual{S_0}\rightarrow B_0$ given by Poincar\'e duality, 
$(\iota|_{S_0})^t:\dual{S_0}\rightarrow\dual{S_0}$ is identified with
$(-1)\circ\iota|_{S_0}: S_0\rightarrow S_0$.

{\it Proof of Claim 1:} We first note that by Lemma 2.2
$(-1)\circ\iota\circ\check\iota:
S_K\rightarrow S_K$ is a holomorphic map, and thus by Torelli for K3 surfaces,
is the identity map if $H((-1)\circ\iota\circ\check\iota):H^2(S_K,\boldz)
\rightarrow H^2(S_K,\boldz)$ is the identity. Thus it is enough
to show that $\check\iota$ and $(-1)\circ\iota$ induce the same
map on cohomology. Thus, from the original construction of $\check\iota$,
we only need to show that $H(-1)=r_P$. But $-1:S_K\rightarrow S_K$
is an involution which clearly leaves $\Pic S_K=P$ invariant, and the
group of anti-invariants is precisely $P^{\perp}$. So $H(-1)|_P=Id_P$
and $H(-1)|_{P^{\perp}}=-Id_{P^{\perp}}$, which is precisely the
definition of $r_P$. $\bullet$

{\it Proof of Claim 2:} Let $P\in\Pone({\bf C})$, $Q=\iota'(P)$, $E_P
=f^{-1}(P)$, $E_Q=f^{-1}(Q)$. The map $\iota|_{E_P}:E_P\rightarrow
E_Q$ is induced by a linear map
$$\dual{(\iota_P^*)}:\dual{H^1(E_P,{\bf R})}\rightarrow
\dual{H^1(E_Q,{\bf R})},$$
which is dual to
$$\iota_P^*:H^1(E_Q,{\bf R})\rightarrow H^1(E_P,{\bf R}),$$
while $\iota^t|_{\dual{E_P}}:\dual{E_P}\rightarrow \dual{E_Q}$ 
is induced by the map
$$\iota_Q^*:H^1(E_P,{\bf R})\rightarrow H^1(E_Q,{\bf R}).$$
Clearly $\iota_P^*\circ \iota_Q^*=Id$, so $\iota_P^*$ and $\iota_Q^*$
are inverses to each other. Furthermore, $\iota|_{E_P}$ and $\iota|_{E_Q}$
are orientation reversing maps, so if $(\ ,\ )_P$ and $(\ ,\ )_Q$ denote
the (skew-symmetric) intersection pairings on $H^1(E_P,{\bf R})$ and
$H^1(E_Q,{\bf R})$ respectively which yield Poincar\'e duality, then
$$(\iota_Q^*\alpha,\iota_Q^*\beta)_Q=-(\alpha,\beta)_P$$
and
$$(\iota_P^*\alpha,\iota_P^*\beta)_P=-(\alpha,\beta)_Q.$$
We can now apply the following lemma to complete the proof.

\proclaim Lemma 2.3. Let $V$ and $W$ be two symplectic vector spaces
of the same dimension, with symplectic forms $\omega_V$ and $\omega_W$
respectively. Let $\phi:V\rightarrow W$ be an isomorphism such that
$\omega_W(\phi(\alpha),\phi(\beta))=-\omega_V(\alpha,\beta)$
for all $\alpha,\beta\in V$. Let $\Psi_V:V\rightarrow \dual{V}$ be
the natural isomorphism given by $\Psi_V(v)(v')=\omega_V(v,v')$,
$\Psi_W:W\rightarrow \dual{W}$ the similarly constructed natural isomorphism.
Then 
$$\Psi_W^{-1}\circ (\phi^{-1})^t\circ \Psi_V=-\phi.$$

Proof:
$$\eqalign{ (
\Psi_W^{-1}\circ (\phi^{-1})^t\circ \Psi_V)(v)&=
(\Psi_W^{-1}\circ (\phi^{-1})^t)(v'\mapsto \omega_V(v,v'))\cr
&=
(\Psi_W^{-1})(w'\mapsto \omega_V(v,\phi^{-1}(w')))\cr
&=
(\Psi_W^{-1})(w'\mapsto -\omega_W(\phi(v),w'))\cr
&=-\phi(v).\quad \bullet\cr}$$ 

This completes the proof of Theorem 2.1. $\bullet$

We would now like to give a more detailed description of the geometry
of the fixed locus of $\iota:S_K\rightarrow S_K$ and how this
changes when we pass to $\check\iota=(-1)\circ\iota$.

For simplicity, we now restrict to the case 
that $\iota$ has $N$ fixed curves, $C_1,\ldots,C_N$, with all
but possibly one being rational. 
(It is easy to check that $Fix(\iota)$ being empty or two elliptic
curves yields a self-mirror K3 and the construction given below still
works.) With respect to the complex structure $K$,
these curves are not holomorphic, however. We give a description
of these curves. We will use heavily the fact that $\iota:S_K\rightarrow
S_K$ is an anti-holomorphic involution, and hence gives a {\it real}
structure on $S_K$; $Fix(\iota)$ is then precisely the
real locus of $S_K$.

First, let $C=Fix(\iota')\subseteq \Pone({\bf C})$ be the fixed locus of
$\iota'$. Now since $Fix(\iota)=C_1\cup\cdots\cup C_N$ is non-empty
and $Fix(\iota)\subseteq f^{-1}(C)$, we must have
$C\not=\phi$, and thus $C$ is the real part of $\Pone({\bf C})$, i.e.
$C\cong S^1$. We can imagine $C$ to be the equator
of $\Pone({\bf C})=S^2$. Thus for each point $P\in C$, $\iota$ maps
$f^{-1}(P)$ to itself, and $\iota:f^{-1}(P)\rightarrow f^{-1}(P)$ is an
anti-holomorphic involution, i.e. $f^{-1}(P)$ has a real structure.
We have
$$Fix(\iota)=\bigcup_{P\in C} Fix(\iota|_{f^{-1}(P)}).$$
Thus we can describe $Fix(\iota)$ as the union of the real parts
of the elliptic curves over $C$. In particular, we have a map 
$Fix(\iota)\rightarrow C$. Now $C'=\sigma\cap f^{-1}(C)$ is fixed by
$\iota$ and $C'$ is also an $S^1$ mapping
isomorphically to $C$ so we have $Fix(\iota)\rightarrow C$ is surjective.
In fact at least one component of $Fix(\iota)$, the one
containing $C'$, maps surjectively to $C$. Let the component of
$Fix(\iota)$ containing $C'$ be $C_1$. Then $C_1$ is not rational,
since it cannot be simply connected. This gives a geometric explanation
of Nikulin's observation that $(S,\iota)$ does not have a mirror
family if $N'=0$ ([3, 17]).

In Figure 1, we depict the map $Fix(\iota)\rightarrow C$ in the case
that $N=3$ and $N'=4$.
\topinsert
$$
\epsfbox{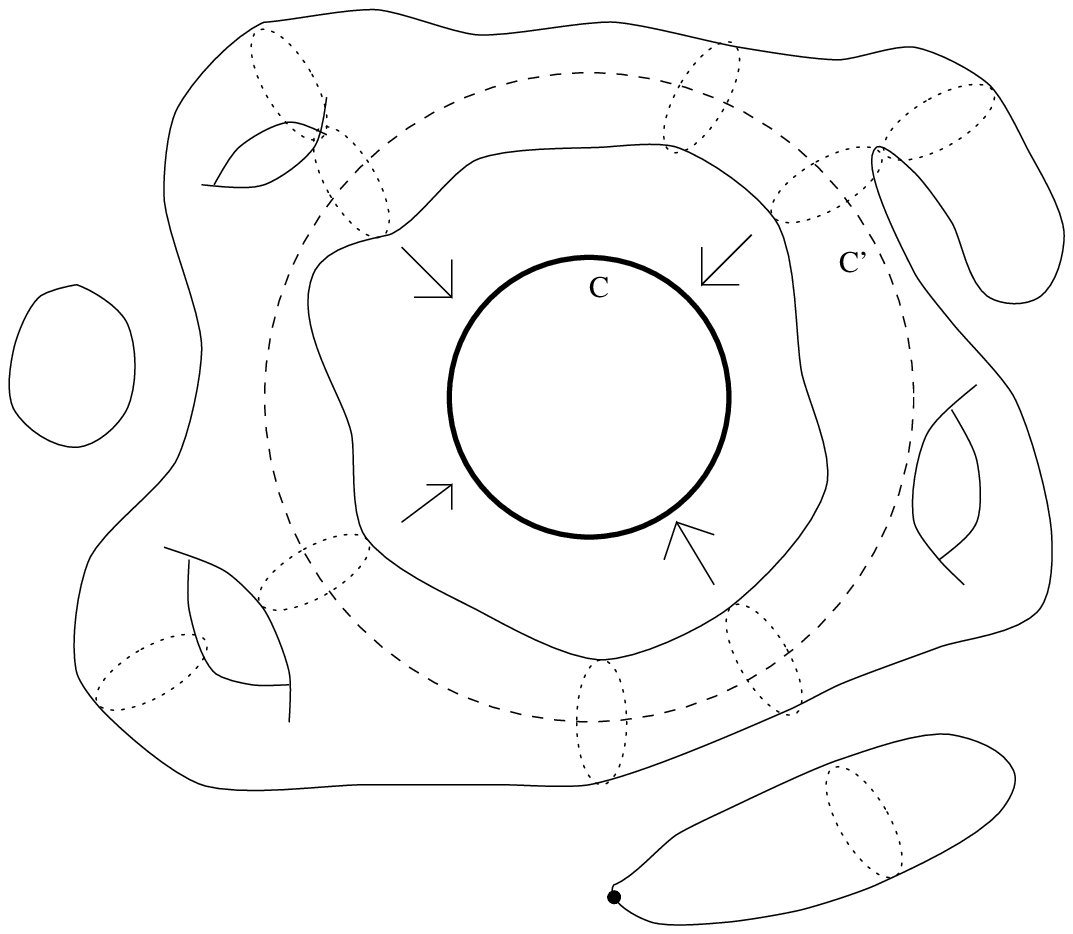}
$$
\bigskip
\centerline{Figure 1}
\endinsert

\topinsert
$$
\epsfbox{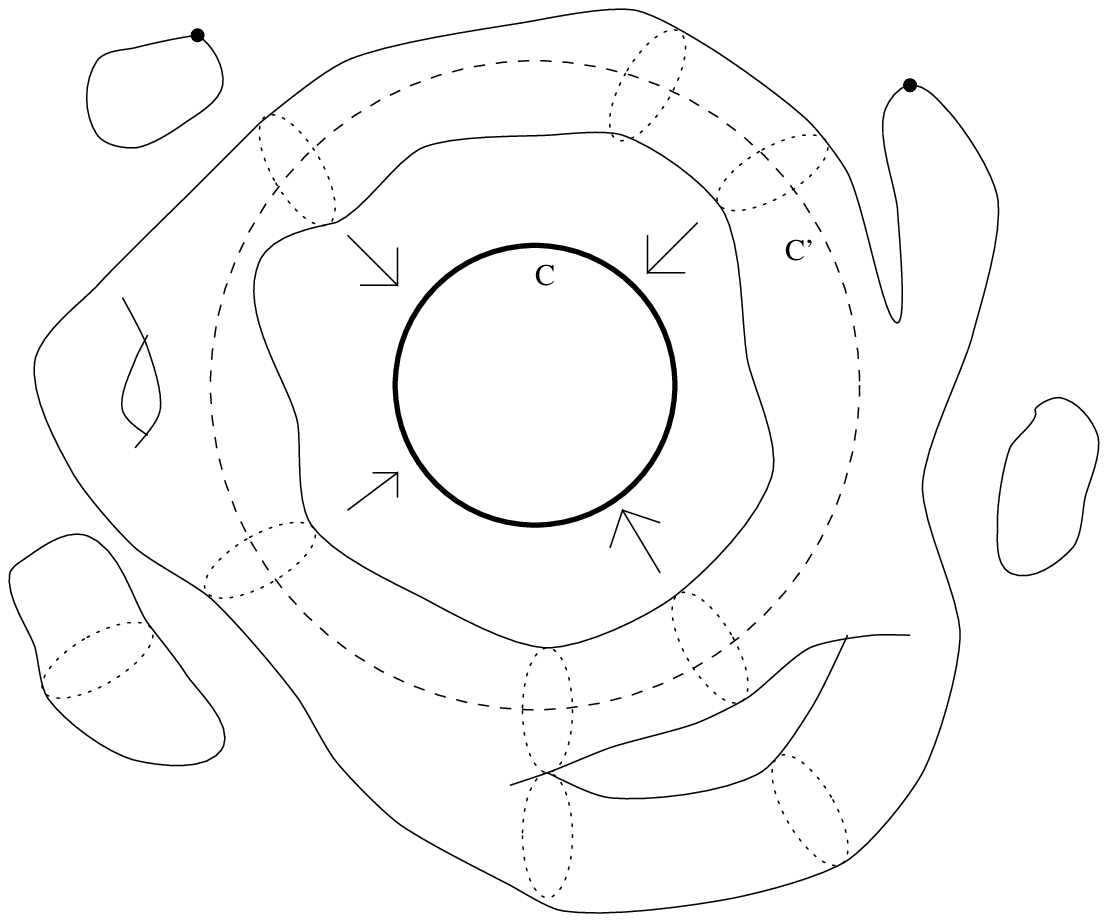}
$$
\bigskip
\centerline{Figure 2.}
\endinsert

We make the following observations to show this picture is reasonably
accurate. First, note that the
real points of a real elliptic curve form either a circle or a union of
two disjoint circles. Thus for each point $P\in B_0\cap C$,
$f^{-1}(P)\cap Fix(\iota)$ consists of one or two circles. This fits
with Figure 1.

If $P\in C$ and $f^{-1}(C)$ is a singular fibre, then from the condition
that $\Pic S_K=U(1)$, $f^{-1}(P)$ must be either of Kodaira
type $I_1$ or type $II$. The
real part of a type $I_1$ fibre is either a figure eight
or a circle plus a point (the node of the $I_1$ fibre). If $f^{-1}(P)$
is a type $II$ fibre, then
the real part of 
$f^{-1}(P)$ in this case must consist topologically of a circle.
Thus if $C^i$ is the subset of $C$ defined by
$$C^i=\{P\in C\ |\ \hbox{$Fix(\iota)\cap f^{-1}(P)$ consists of $i$ disjoint
circles}\},$$
$i=1,2$, then for any $P\in C-C^1\cup C^2$, $f^{-1}(P)$ is of type $I_1$,
$Fix(\iota)\cap f^{-1}(P)$ consists either of a figure eight or
a circle plus point, and $P$ is a boundary point between $C^1$ and $C^2$.
This fits the picture in Figure 1, where we draw various fibres of
$Fix(\iota)\rightarrow C$ as dotted circles.

Now we have shown that $\check\iota=(-1)\circ\iota$. Note that
if $E_0$ is a nodal cubic with anti-holomorphic involution $\iota:E_0\rightarrow
E_0$, then $-\iota:E_0\rightarrow E_0$ 
is also anti-holomorphic, and if the fixed
locus of $\iota$ was a figure eight, it is easy to see that the fixed
locus of $-\iota$ is a circle and a point. One way of seeing
this is as follows. Suppose we have an elliptic curve $E$ with
periods $\{1,\alpha i\}\subseteq {\bf C}$ with anti-holomorphic involution 
$\iota$ given by complex conjugation, $z\mapsto \bar z$. Then $Fix(\iota)$
is given by the two circles $\im z=0$ and $\im z=\alpha /2$. If we let
$\alpha$ go off to $\infty$ to obtain $E_0$, then this is the same as
shrinking the cycle given by $\im z=0$ to a point. Thus the fixed
part of $E_0$ under $\iota$ consists of a circle and a point. If we now
consider the involution $-\iota$ on $E$, it has fixed locus consisting
of the circles $\re z=0$
and $\re z=1/2$. Shrinking $\im z=0$ now yields a figure eight as fixed locus
of $\iota$ on $E$. In fact, one can show that the elliptic curves
over $C^2\subseteq C$ will indeed be of this form, with purely imaginary
period.

Thus, passing between
$\iota$ and $\check \iota$ has the effect of interchanging figure eights
with circles plus points. Carrying out this transformation in Figure
1 by replacing each figure eight with a circle plus point, and each
circle plus point with a figure eight, we obtain Figure 2. We 
break the ``handles'' off from $C_1$ and ``sew'' the spheres
$C_2,\ldots,C_{N}$ into $C_1$. This has the effect of interchanging
the numbers $N$ and $N'$, as was proved by other methods [13, 3, 17].

\bigskip
{\hd \S 3. Borcea-Voisin examples of mirror Calabi-Yau threefolds}

A construction of mirror symmetry for a particular class of Calabi-Yau
threefolds was  proposed by Borcea [3] and Voisin [17].
The threefolds in question are obtained from K3 surfaces
$S$ with an involution $\iota$ (acting by negation on the holomorphic
2-form) as described in \S 2.
For simplicity, we shall neglect the possibility of the fixed locus
of $\iota$ being empty or
two elliptic curves, where it may be checked that both the
K3 surface and also the Calabi-Yau threefold constructed below are self-mirror.
Thus we assume that $\iota$ has $N$ fixed curves, $N-1$ of  which are rational
and one being of genus $N'\ge 0$.
As noted in \S 2, if $N>0$, then $N'>0$ is
a necessary condition for there to be
a mirror of $(S,\iota)$, and hence for the Borcea-Voisin construction to work.

If we now take any elliptic curve $A$ with $j$ denoting its standard involution
given by negation, the action of $(\iota,j)$ on $S\times A$ has $4N$
fixed curves corresponding to the fixed curves of $\iota$ and the
$2$-torsion points of $A$. The quotient $Y=S\times A/(\iota,j)$
therefore has $4N$ curves of $A_1$ singularities, each of which
may be resolved by a single blowing-up. It is clear that this
resolution $X$ is a Calabi-Yau threefold. A calculation of the Hodge
numbers [3, 17] shows that
$$\eqalign{h^{1,1}(X)&=11+5N-N'\cr
h^{2,1}(X)&=11+5N'-N.\cr}$$
We now assume that $(S,\iota)$ has a mirror $(\check S,\check\iota)$, i.e.
that the transcendental lattice $T$ of $S$ contains a hyperbolic plane
$P$. Necessary and sufficient conditions for this are given in [17], (2.5).
As was noted in \S 2, the fixed locus of $\check \iota$
consists of $N'$ smooth curves, one of genus $N$ and the others all
rational.

The mirror Calabi-Yau threefold is constructed in the obvious way: we have
the involution $(\check \iota,j)$ acting on $\check S\times A$
(recalling that an elliptic curve is self-mirror) and set $\check X$ to
be the minimal resolution of $\check S\times A/(\check\iota,j)$. It is then
clear that $h^{1,1}(X)=h^{2,1}(\check X)$, and $h^{2,1}(X)=h^{1,1}(\check
X)$ from the formulae given above, and that the
Euler characteristics are related by $e(X)=12(N-N')=-e(X')$. Further
properties of these mirror pairs are checked in [17], although it is
observed that the construction in fact yields only a slice of the mirror
map, as $X$ and $X'$ are restricted in moduli to those arising as quotients
in the way described above. We'll return to this point later
when discussing our construction of the mirror.

In this section, we verify that the Borcea-Voisin construction is consistent
with the recipe suggested in [16], at least if we allow a degenerate
K\"ahler-Einstein metric on $X$, corresponding to degenerating the
K\"ahler class to the wall of the K\"ahler cone consisting of the pullbacks
of ample classes on $Y$. More specifically, the choice of such a nef 
class on $X$ determines a K\"ahler class on $S\times A$, which in turn
determines K\"ahler classes on $S$ and $A$. These in turn determine 
a Ricci-flat metric on $S$ and a flat metric on $A$. Since the K\"ahler
class on $S\times A$ is invariant under our involution, so too is the
K\"ahler-Einstein metric we have constructed, which therefore descends to
give an orbifold metric on $Y$ and a degenerate metric on
$X$, which is a K\"ahler-Einstein metric on the complement of the exceptional
locus.

{\it Remark.} By degenerating the K\"ahler class to a pullback of an ample
class on $Y$, we are imposing the condition that the nef class
lies in a codimension $4N$ wall of the K\"ahler cone--- we saw above that there
are $4N$ exceptional divisors on $X$. If we also restrict the class $B$
to come from pulling back a class on $Y$, we are imposing a complex 
codimension $4N$ condition on the degeneration of the complexified
K\"ahler class. On the mirror, we obtain from the Borcea-Voisin construction
a Calabi-Yau threefold $\check X$ containing four ruled surfaces over a curve
of genus $N$; this is well-known to impose a codimension $4N$ condition
on the complex moduli space. By restricting to these subsets of complex
and K\"ahler moduli spaces, all calculations may essentially be reduced
to ones on $S\times A$ (cf \S 4 where we describe the mirror map); 
this is the sense in which we are only seeing
a slice of the mirror map, but one which is sufficient to determine
the mirror as a topological manifold.

We will begin by applying the recipe of [16] with the degenerate 
K\"ahler-Einstein metric defined above. By extension of previous terminology,
we shall however say that a torus fibration on $X$ 
is special Lagrangian (with respect to this degenerate metric) if the
general fibre is a $T^3$ contained in the complement of the exceptional
locus, and is special Lagrangian with respect to the K\"ahler-Einstein
metric on this complement.

Choosing a K\"ahler class on $S$, we have the fibration $S\rightarrow S^2$
(corresponding to the elliptic fibration $S_K\rightarrow\Pone$ when we
change the complex structure) whose general fibre is special Lagrangian
with respect to the K\"ahler-Einstein metric on $S$. Similarly, we can choose 
a special Lagrangian $T^1$-fibration on $A$ with respect to a flat metric;
if for instance $A$ is taken with periods $\{1,\tau\}\subseteq {\bf C}$
with coordinate $z$,
then
the curves $\im z=constant$ yield an appropriate fibration. Thus with
respect to the K\"ahler-Einstein metric on $S\times A$, we obtain a special
Lagrangian torus fibration $\tilde\pi:S\times A\rightarrow S^2\times S^1$
(whose fibres are special Lagrangian 3-tori except over a finite number of
copies of
$S^1$, corresponding to the points of $S^2$ over which the K3 surface
has singular fibres). As noted in \S 2, 
for
a general choice of the period ${\bf C}\Omega$ in $D_M'\backslash \Delta_M$
and K\"ahler class $\omega$, by construction we have
$\Pic(S_K)=P$ 
the hyperbolic plane. So the elliptic fibration $S_K\rightarrow\Pone$
has only singular fibres of Type $I_1$ and $II$,
and a standard Euler characteristic computation shows that
$$\#\{\hbox{Type $I_1$ fibres}\}+2\cdot\#\{\hbox{Type $II$ fibres}\}=24.$$

Passing to the quotient $Y$, we obtain a special Lagrangian torus
fibration with respect to the orbifold metric on $Y$ corresponding
to the K\"ahler-Einstein metric on $S\times A$. The base of this
fibration is $B=(S^2\times S^1)/(\iota',j)$. Here we recall that the action
$\iota'$ on the base $S^2$ was induced via the action of $\iota$ on the section
$\sigma$ and was anti-holomorphic. It may therefore be taken as reflection
about the equator in $S^2$. The action of $j$ on the base $S^1$ of the
special Lagrangian fibration on $A$ 
may be taken to be reflection
in the real axis, thinking of $S^1$ as embedded in the complex plane
in the natural way. The general fibre of $\bar\pi:Y\rightarrow B$
is still a special Lagrangian $T^3$ (with respect to the
orbifold metric on $Y$). The induced map $\pi:X\rightarrow B$
is therefore a special Lagrangian torus fibration
(with respect to the degenerate metric on $X$ described above).
We shall describe $B$ and the discriminant locus of $\pi:X\rightarrow
B$ in more detail below. However, we first note that the
Borcea-Voisin construction corresponds precisely with the recipe proposed in
[16].

\proclaim Proposition 3.1. The mirror $\check X$ of $X$ as constructed
above is a compactification of the dual torus fibration of
$p:X_0\rightarrow B_0$, where $B_0$ is the complement of the discriminant
locus of $\pi$ in $B$, $X_0=\pi^{-1}(B_0)$, $p=\pi|_{X_0}$ and
the dual fibration is defined as before to be $R^1p_*{\bf R}/R^1p_*\boldz$, as
a $T^3$ fibration over $B_0$. 

Proof: This follows from the results of \S 2. We have a commutative
diagram
$$\matrix{S\times A&\mapright{}&Y\cr
\mapdown{\tilde\pi}&&\mapdown{\bar\pi}\cr
S^2\times S^1&\mapright{}&B\cr}$$
We saw in \S 2 that dualizing $(S,\iota)$ as a torus fibration over
$S^2$ yielded the mirror $(\check S,\check \iota)$, with $\check\iota$
the transpose of $\iota$. Dualizing the $T^1$
fibration of $A$ over $S^1$ just recovers $A$ again, and 
the transpose of $j$ is again $j$. Thus, taking the dual
of the torus fibration $S\times A\rightarrow S^2\times S^1$, we recover
$\check S\times A$ with involution $(\check\iota,j)$. Therefore
dualizing and taking the quotient by $(\check\iota,j)$ will be
the same (at least over $B_0$) as taking the dual of the smooth torus
fibration on $X_0$. Thus $\check X$ is a smooth compactification of
this dual torus fibration. $\bullet$

{\it Remarks.} (1) The construction ensures that the smooth fibres of 
the dual fibration are also special Lagrangian with respect to 
the (degenerate) metric on $\check X$ deduced from the relevant
K\"ahler metrics
on $\check S$ and $A$.

(2) We conjecture that in general, when faced with a special Lagrangian
torus fibration on one Calabi-Yau, the (allowable) compactification
of the dual should be determined uniquely by knowledge of the monodromies.

(3) The recipe in [16] also suggests a method for putting a complex
structure on the smooth part of the dual fibration--- we shall discuss
this further for the Borcea-Voisin examples in \S 4.

For the remainder of this section, we concentrate on the topology of the 
manifolds $X$ and $\check X$, viewed as fibre spaces over $B$ with general 
fibres being $3$-tori.

\proclaim Proposition 3.2. The base $B$ may be identified topologically
as $S^3$.

Proof: Recall that the induced involution on $S^2\times S^1$ consists
of reflection in the equator on $S^2$ and reflection in
the real axis on $S^1$. This clearly has fixed locus $S^1\times S^0$.
To see that topologically the quotient is $S^3$, we use real coordinates,
writing $S^2\subseteq{\bf R}^3$ as $x^2+y^2+z^2=1$ and $S^1\subseteq{\bf R}^2$
as $u^2+v^2=1$. We take the involution to be the one changing the signs of
both the $z$ coordinate and the $v$ coordinate. The invariant polynomials
are therefore $X=x$, $Y=y$, $Z=z^2$, $U=u$, $V=v^2$, $W=zv$. The quotient
can therefore be realised as the subset of ${\bf R}^6$ with equations
$$\eqalign{X^2+Y^2+Z&=1\quad (Z\ge 0)\cr
U^2+V\hfil&=1\quad (V\ge 0)\cr
W^2\hfil &=ZV.\cr}$$
Eliminating $V$ and $Z$, $B$ is identified as the subset of ${\bf R}^4$
given by 
$$(1-U^2)(X^2+Y^2)+U^2+W^2=1$$
where $U^2\le 1$ and $X^2+Y^2\le 1$. Setting $X=R\cos\theta$, $Y=
R\sin\theta$ ($-1\le R\le 1, 0\le\theta<\pi$), we obtain an equation
$W^2=(1-U^2)(1-R^2),$ where $-1\le U\le 1$ and $-1\le R\le1$, i.e.
a double cover of the square. This double cover has singularities
precisely at the corners of the square (locally $W^2=ST$
with $S\ge 0$, $T\ge 0$),
but is topologically just an $S^2$. The quotient we seek is therefore
obtained by rotating this $S^2$ about $R=0$, thus obtaining an $S^3$
topologically. $\bullet$

We can now identify the discriminant locus on $B$ of our torus fibration
$X\rightarrow B$. One part consists of the image of
the fixed locus of the involution $(\iota',j)$ acting on $S^2\times S^1$;
this we saw consisted of two disjoint
copies of $S^1$, and more explicitly was the
locus on $B$ swept out by the singular points ($R=\pm 1, U=\pm 1, W=0$)
when we rotate about $R=0$. The second part of the discriminant locus
occurs because of singular fibres of the special Lagrangian torus
fibration $S\rightarrow S^2$ on the K3. For each $P\in S^2$ corresponding
to a singular fibre, we have a component of the discriminant locus on $B$,
namely the image of $P\times S^1$ under the quotient map.
If $P$ is not a fixed point of the induced involution on $S^2$, this is just
an $S^1$ on $B$ disjoint from the two copies of $S^1$ 
previously identified. If however
$P$ is a fixed point, then the image of $P\times S^1$ will be an interval
whose endpoints lie on these two circles, i.e. the two
components of the image of the fixed locus
$S^1\times S^0\subseteq S^2\times S^1$. Schematically, we can now picture
the discriminant locus on $S^3$ consisting of two disjoint circles,
$\Gamma_1$ and $\Gamma_2$, corresponding to the fixed locus of
the involution on $S^2\times S^1$, and one extra component (either an
interval or a circle) for each singular fibre of $S\rightarrow S^2$.
This is pictured in Figure 3.
\topinsert
$$
\epsfbox{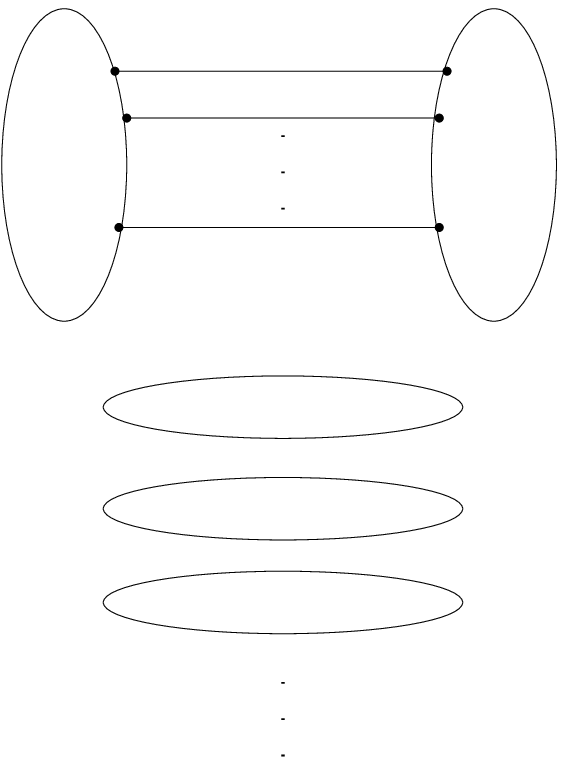}
$$
\bigskip
\centerline{Figure 3.}
\bigskip
\endinsert

The fibre of $Y\rightarrow B$ above a general point of one of the
two circles $\Gamma_i$ is a 3-torus $T$ modulo involution and when passing to $X$
we replace in this fibre 2 or 4 copies of $S^1$ by $S^2$ bundles over $S^1$.
The fibres therefore do not make any contribution to the Euler
characteristic $e(X)$. The fibre of $Y\rightarrow B$ over
a general point of a component corresponding to a singular fibre
$Z$ of $S\rightarrow S^2$ will be $Z\times S^1$, and therefore also
does not contribute to $e(X)$. The Euler characteristic may therefore
be explained purely in terms of the points of intersection in the above
picture, i.e. the fibres of $Y\rightarrow B$ 
of the form $Z\times S^1/(\iota,j)$ 
with $Z$
a fixed singular fibre on $S\rightarrow S^2$. If $Z_0$ denotes the
fixed locus of the induced involution on $Z$, we know that $Z_0$ is
either a figure eight or a circle plus a point (for $Z$ a type $I_1$
fibre) or a singular circle (for $Z$ a type $II$ fibre). The corresponding
fibre of $X$ is therefore obtained by replacing
$Z_0\times S^0$ in $Z\times S^1/(\iota,j)$ 
by an $S^2$-bundle over $Z_0\times S^0$.
Thus when $Z$ is of type $II$, the local contribution is still trivial.

Let us suppose therefore that $Z$ is of type $I_1$ and that $Z_0$ is a
figure eight. The fixed locus of $\iota$ on $S\times A$ therefore contains
four copies of $Z_0$ (one for each 2-torsion point of $A$), two copies above
each of the two intersection points on the discriminant locus. Thus
for each of the two fibres of the form $Z\times S^1$, we are removing
two copies of $Z_0$ in the manifold $S\times A$, thus increasing the
Euler characteristic by 2 for each fibre. Therefore on the quotient $Y$,
the Euler characteristic has been increased by 1 for each fibre. When we
pass to $X$, we replace each copy of $Z_0$ by an $S^2$-bundle over
$Z_0$, thus decreasing the Euler characteristic by $4$ (for each
fibre). Over each point of intersection, the fibre of $X$ over $B$ therefore
has a net contribution of $-3$ to $e(X)$, and so the singular fixed fibre
$Z$ of $S\rightarrow S^2$ accounts for a total contribution of $-6$. A similar
calculation shows that when $Z_0$ is a circle plus a point, the net contribution
from each of the two fibres $Z\times S^1/(\iota,j)$ will be $+3$, making a total
contribution of $+6$ to $e(X)$.

We saw however in \S 2 that there will be $2(N-1)$ fibres $Z$ of $S\rightarrow
S^2$ where $Z_0$ is a circle and a point and corresponding
to the $N-1$ fixed rational curves on $S$, and $2(N'-1)$ fibres
$Z$ where $Z_0$ is a figure eight and corresponding to holes of
the genus $N'$ fixed curve (excluding the hole in the middle). There
may also be other fibres $Z$ with $Z_0$ one of the above two types,
but the two types will have to occur in equal numbers. Thus the net
contribution we obtain to $e(X)$ from the singular fibres is
$2(N-1)\cdot 6-2(N'-1)\cdot 6=12(N-N')$, in agreement with the 
value of $e(X)$ previously calculated. We recall also however from
\S 2 that when dualizing the two possibilities for $Z_0$ are switched
around--- recall that dualizing changes the original involution $\iota$
to $\check\iota$. This confirms immediately that the effect of dualizing
(i.e. passing from $X$ to $\check X$) on the Euler characteristic is merely
to change its sign.
\vfill
\eject
{\hd \S 4. Mirror Maps and Conjectures.}

So far, we have not dealt with the most mysterious aspect of
the SYZ construction: namely, how does one place a complex
structure on the dual of a special Lagrangian torus fibration?
This is one of the key questions raised in [16].
Section 3 of [16]
gives a local construction of an almost complex structure, which is then
shown
to be integrable. An alternative approach is explained
in [9]. This will not however (even locally) define the correct complex
structure on the mirror, since the physics predicts that there
will also be instanton corrections coming from the singular
fibres of the fibration [16, \S 2].  

We give one approach, which works in the K3 and Borcea-Voisin 
examples, and we expect it to work more generally. 
In brief, given a fixed complex structure and
K\"ahler form on a Calabi-Yau threefold $X$ which yields a special Lagrangian
$T^3$-fibration $f:X\rightarrow B$ with special Lagrangian section, 
we then conjecture that the
mirror map can be expressed in terms of the Leray spectral
sequence arising from $f$. For a fixed such $f$, one then obtains a map
between the K\"ahler moduli of $X$ and the complex moduli of $\check X$.
This map should not depend on the initial choice of $f$, however.
We will formulate this in a more precise conjecture at the
end of this section; we first study explicit examples to give
motivation for the conjecture.

Let $S$ be a K3 surface as in \S 1 with $P\cong
U(1)\subseteq T$, 
and
choose complex and K\"ahler structures as in \S 1 so that 
there is a corresponding special Lagrangian 2-torus fibration
$f : S \to S^2$ with special Lagrangian section $\sigma$ with
$E' = E+ \sigma $ (see Remark following (1.3)). 
We use the same letters to represent the cohomology
classes of these curves.
We can assume in what follows that $f:S\rightarrow S^2$ has
no reducible fibres by choosing sufficiently general $\omega$ and $\Omega$
on $S$.

\proclaim Proposition 4.1. The Leray spectral sequence
$$H^i(S^2,R^jf_*{\bf Q})\Rightarrow H^n(S,{\bf Q})$$
degenerates at the $E_2$ term, and the non-zero terms are shown in
the following diagram:
$$\matrix{{\bf Q}&0&{\bf Q}\cr
0&H^1(S^2,R^1f_*{\bf Q})&0\cr
{\bf Q}&0&{\bf Q}\cr}$$

Proof: This all follows from [5], \S 1, (1.4) and following,
where ${\bf Q}=H^0(S^2,R^2f_*{\bf Q})$
is generated by $\sigma$ (or equivalently by $E'$) and ${\bf Q}
=H^2(S^2,f_*{\bf Q})$ is generated by $E$. $\bullet$

This yields a filtration
$$0\subseteq F_0\subseteq F_1\subseteq F_2=H^2(S,{\bf Q})$$
such that
$$\eqalign{
F_0&\cong H^2(S^2,f_*{\bf Q})\cong {\bf Q}E\cr
F_1/F_0&\cong H^1(S^2,R^1f_*{\bf Q})\cr
F_2/F_1&\cong H^0(S^2,R^2f_*{\bf Q})\cong {\bf Q}E'\cr}$$

This filtration allows us to give a very natural alternative construction
for  the mirror maps
$\check\phi:T_M\rightarrow D_{\check M}$ and
$\check\phi^{-1}:D_{\check M}\rightarrow T_M$
as follows. 
First we can identify $S \to S^2$ as a 
compactification also of the dual fibration via Poincar\'e 
duality as in \S 1, and thus ignore the dualizing.
Now consider ${\bf C}\check\Omega\in D_{\check M}$. 
After an appropriate scaling, $\im\check\Omega\in
F_1\otimes_{\bf Q}{\bf R}$ and $\re\check\Omega\in E'+F_1\otimes_{\bf Q}
{\bf R}$.
Thus $\check\Omega-E'\in F_1\otimes_{\bf Q}{\bf C}$, and we take
$\check\phi^{-1}(\check\Omega)=\check\Omega-E'\in
(F_1/F_0)\otimes_{\bf Q}{\bf C}$. As elements of $(F_1/F_0)\otimes_{\bf Q}
{\bf C}$, this will coincide with $\check\phi^{-1}(\check\Omega)=
B+i\omega$ as defined in \S 1.

Thus, conversely, to define $\check \phi$, we take the class of $B+i\omega
\in T_M$ in $F_1/F_0\otimes_{\bf Q}{\bf C}$ and lift this to a class
$\alpha$ in $F_1\otimes_{\bf Q}{\bf C}$ in such a way that 
$\alpha+E'\in D_{\check M}$. We claim there is a unique such
lifting. This follows immediately from the proof of Proposition 1.1.
We then define $\check\phi(B+i\omega)=\alpha+E'$.

{\it Remark.} A similar construction produces the mirror map
for elliptic curves, and more generally, for complex tori of
any dimension. We sketch this for elliptic curves, as we will
need it below. If $g:A\rightarrow S^1$ is a special Lagrangian
fibration with respect to some fixed complex and K\"ahler structures, 
choose real coordinates $x$ and $y$ on $A$ so that
$A\cong {\bf R}^2/((1,0)\boldz+(0,1)\boldz)$ and $g:A\rightarrow S^1$
is given by $(x,y)\mapsto y$. The Leray spectral sequence for $g$
at the $E_2$ level is
\def\Q{{\bf Q}}
$$\matrix{H^0(S^1,R^1g_*\Q)&H^1(S^1,R^1g_*\Q)\cr
H^0(S^1,g_*\Q)&H^1(S^1,g_*\Q)\cr}$$
which clearly degenerates. We can think of $H^0(S^1,R^1g_*\Q)$ as 
being generated by $s_x=\{x=constant\}$, a section of $g$, and
$H^1(S^1,g_*\Q)$ is generated by $s_y=\{y=constant\}$, a fibre of
$g$. The period domain for $A$ is 
\def\C{{\bf C}}
$$D_A=\{\C\Omega\in\P(H^1(A,\C))\ |\ i\Omega\wedge\bar\Omega>0\}.$$
By normalizing $\Omega$ for $\C\Omega\in D_A$, we can always write
$\Omega$ uniquely as $s_x+\tau s_y$ for $\tau\in\H$, the upper half plane.
This gives an isomorphism between $D_A$ and the tube domain
$$T_A=\{B+i\omega\in H^2(A,\C)=\C\ |\ \omega>0\}$$
via $\phi_A:T_A\rightarrow D_A$ defined by 
$$\phi_A(B+i\omega)=s_x+(B+i\omega)s_y.$$
The act of dualizing the $S^1$ fibration $g:A\rightarrow S^1$ to
obtain $\dual{g}:\dual{A}\rightarrow S^1$ has the effect of interchanging
$R^1g_*\Q$ and $g_*\Q$, i.e. interchanging the two rows of the spectral
sequence. Thus, given a complexified K\"ahler class $B+i\omega\in T_A
\subseteq H^2(A,\C)=H^1(S^1,R^1g_*\C)$, we obtain an element of
$H^1(S^1,\dual{g}_*\C)\subseteq H^1(\dual{A},\C)$. Adding $s_x$ to this
element yields $\phi_A(B+i\omega)$. This is precisely the same
recipe as was performed above in the K3 case.

We will show a similar construction works for the case of Borcea-Voisin
threefolds, and we believe that this construction is similar
in spirit to the ideas for putting complex structures on the dual fibration
given in [16].

We now study the fibration $\bar\pi:Y\rightarrow S^3$ constructed
in \S 3, given a choice of $(S,\iota)$ of K3 surface with involution
and elliptic curve $A$, with fixed complex and K\"ahler structures.

\proclaim Proposition 4.2. The Leray spectral sequence
$$H^i(S^3,R^j\bar\pi_*{\bf Q})\Rightarrow H^n(Y,{\bf Q})$$
degenerates at $E_2$, and the non-zero terms of $E_2$ are shown in the
following diagram:
$$\matrix{{\bf Q}&0&0&{\bf Q}\cr
0&\check M_{\bf Q}\oplus {\bf Q}&
M_{\bf Q}\oplus {\bf Q}&0\cr
0&M_{\bf Q}\oplus {\bf Q}&
\check M_{\bf Q}\oplus {\bf Q}&0\cr
{\bf Q}&0&0&{\bf Q}\cr}$$
where $M_{\Q}=M\otimes_{\boldz} \Q$ and $\check M_{\Q}=\check M\otimes_{\boldz}
\Q$. 

Proof. First we compute $H^i(Y,\Q)$.  Let $G=\boldz/2\boldz$ be
the group acting on $S\times A$ generated by the involution
$(\iota,j)$. Then 
$$H^r(Y,\Q)=H^r(S\times A,\Q)^G=\left(
\bigoplus_{i+j=r} H^i(S,\Q)\otimes H^j(A,\Q)\right)^G$$
by K\"unneth. Thus
$$\eqalign{
H^0(Y,\Q)&=\Q\cr
H^1(Y,\Q)&=0\cr
H^2(Y,\Q)&=H^2(S,\Q)^+\otimes H^0(A,\Q)
\oplus H^0(S,\Q)\otimes H^2(A,\Q)= M_{\Q}\oplus \Q\cr
H^3(Y,\Q)&=H^2(S,\Q)^-\otimes H^1(A,\Q)\cr
H^4(Y,\Q)&=H^2(S,\Q)^+\otimes H^2(A,\Q)\oplus H^4(S,\Q)\otimes H^0(A,\Q)
=M_{\Q}\oplus \Q\cr
H^5(Y,\Q)&=0\cr
H^6(Y,\Q)&=\Q.\cr}$$
\def\bt{\boxtimes}

Now let $f:S\rightarrow S^2$ be the 2-torus fibration as before and
$g:A\rightarrow S^1$ the chosen circle fibration on $A$, as in the remark
above, so that $\tilde\pi:S\times A\rightarrow S^2\times S^1$
is $f\times g$. Denote, for sheaves $\F$ and $\G$ on $S^2$ and $S^1$
respectively, $\F\bt\G:= p_1^*\F\otimes p_2^*\G$ with $p_1$ and $p_2$
the first and second projections of $S^2\times S^1$ onto $S^2$ and
$S^1$ respectively. Then again by K\"unneth,
$$\eqalign{
\tilde\pi_*\Q&= f_*\Q\bt g_*\Q\cr
R^1\tilde\pi_*\Q&=
R^1f_*\Q\bt g_*\Q\oplus f_*\Q\bt R^1g_*\Q\cr
R^2\tilde\pi_*\Q&=
R^2f_*\Q\bt g_*\Q\oplus R^1f_*\Q\bt R^1g_*\Q\cr
R^3\tilde\pi_*\Q&=
R^2f_*\Q\bt R^1g_*\Q.\cr}$$
Now
$$H^i(S^3,R^j\bar\pi_*\Q)=(H^i(S^2\times S^1,R^j\tilde\pi_*\Q))^G.$$
From this, the values of $H^i(S^3,R^j\bar\pi_*\Q)$ for $j=0$
and $3$ given in Proposition 4.2 are clear, and we check the middle two
rows.

We have
$$\eqalign{
H^0(S^3,R^1\bar\pi_*\Q)&=(H^0(S^2\times S^1,R^1\tilde\pi_*\Q))^G\cr
&=(H^0(S^2,R^1f_*\Q)\otimes H^0(S^1,g_*\Q)
\oplus H^0(S^2,f_*\Q)\otimes H^0(S^1,R^1g_*\Q))^G.\cr}$$
Note that $H^0(S^2,R^1f_*\Q)=0$ by Proposition 4.1, and
$G$ acts on 
$H^0(S^2,f_*\Q)\otimes H^0(S^1,R^1g_*\Q)$ trivially on the first factor
and by negation on the second factor, as is clear by inspecting the
Leray spectral sequence for $g$, so we obtain
$H^0(S^3,R^1\bar\pi_*\Q)=0$.

Next
$$\eqalign{H^1(S^3,R^1\bar\pi_*\Q)=&(H^1(S^2\times S^1,R^1\tilde\pi_*\Q))^G\cr
=&(H^0(S^2,R^1f_*\Q)\otimes H^1(S^1,g_*\Q)\cr
&\oplus H^1(S^2,R^1f_*\Q)\otimes H^0(S^1,g_*\Q)\cr
&\oplus H^0(S^2,f_*\Q)\otimes H^1(S^1,R^1g_*\Q)\cr
&\oplus H^1(S^2,f_*\Q)\otimes H^0(S^1,R^1g_*\Q))^G\cr
=& H^1(S^2,R^1f_*\Q)^+\otimes H^0(S^1,g_*\Q)
\oplus H^0(S^2,f_*\Q)\otimes H^1(S^1,R^1g_*\Q)\cr
=& M_{\Q}\oplus \Q\cr}$$
Continuing as above, we fill in the rest of the table, noting the
particular isomorphisms
$$\eqalign{E^{2,1}_2\cong &H^1(S^2,R^1f_*\Q)^-\otimes H^1(S^1,g_*\Q)
\oplus H^2(S^2,f_*\Q)\otimes H^0(S^1,R^1g_*\Q)\cr
\cong &\check M_{\Q}\otimes \Q s_y\oplus \Q E\otimes \Q s_x\cr}$$
and
$$\eqalign{E^{1,2}_2\cong &
H^0(S^2,R^2f_*\Q)\otimes H^1(S^1,g_*\Q)
\oplus 
H^1(S^2,R^1f_*\Q)^-\otimes H^0(S^1,R^1g_*\Q)
\cr
\cong &\Q E'\otimes \Q s_y\oplus \check M_{\Q}\otimes \Q s_x,\cr}$$
as well as
$$\eqalign{E^{0,3}_2\cong& H^0(S^2,R^2f_*\Q)\otimes H^0(S^1,R^1g_*\Q)\cr
\cong & \Q E'\otimes \Q s_x.\cr}$$

This spectral sequence must then degenerate at the $E_2$ term, giving the
known values of $H^n(Y,\Q)$ computed above. In fact, this could have
been deduced from the fact that the Leray spectral sequences for both
$f$ and $g$ degenerate. $\bullet$

{\it Remark.} One can also compute the spectral sequence for $\pi:X\rightarrow
B$, and one finds one needs to modify the $E_2^{0,2}$, $E_2^{1,2}$,
$E_2^{0,3}$ and $E_2^{1,3}$ terms appropriately. However, for our
purposes, it appears that if one is using a special Lagrangian 3-torus
fibration associated to a degenerate K\"ahler class, it is more appropriate
to work with the singular Calabi-Yau threefold.

We can now use this spectral sequence as before to produce the mirror
map for the Borcea-Voisin threefolds. We first define suitable
tube and period domains. As mentioned in \S 3, we will only be constructing
the mirror maps on certain slices of complex and K\"ahler moduli. We will
restrict attention to those complex moduli of $X$ which can be described
as a resolution of a quotient of $S\times A$, and those K\"ahler
forms which come from pull-backs of K\"ahler forms on
$Y$. Recalling that 
$$H^3(Y,\C)=
H^2(S,\C)^-\oplus H^1(A,\C)$$
and
$$H^2(Y,\C)=
H^2(S,\C)^+\oplus H^2(A,\C),$$
we define
$$\eqalign{D_{Y}=\{\C\Omega\in \P(H^3(Y,\C))
\ |&\ \Omega=\Omega_S\otimes \Omega_A\in H^2(S,\C)^-\otimes H^1(A,\C)\cr
&\hbox{with $\C\Omega_S\in D_M$ and $\C\Omega_A\in D_A$.}\}\cr}$$
and
$$\eqalign{T_{Y}=\{B+i\omega\in H^2(Y,\C)
\ | &\ B+i\omega=(B_1+i\omega_1,B_2+i \omega_2)\in H^2(S,\C)^+\oplus H^2(A,\C),
\cr
& B_1+i\omega_1\in T_M, B_2+i\omega_2\in T_A\}.\cr}$$
We note that by Torelli on $S$ and $A$, the period domain $D_Y$ parametrizes
complex structures on $S\times A$ which are invariant under $G$, and hence
parametrizes
complex structures on $Y$.

The Leray spectral sequence associated to the given fibration
$\bar\pi:Y\rightarrow S^3$ yields a filtration
$$0\subseteq F_0\subseteq F_1\subseteq F_2\subseteq F_3=H^3(Y,\Q)$$
with $F_i/F_{i-1}\cong E_2^{3-i,i}$. The values of these quotients are
given explicitly in the proof of Proposition 4.2. The mirror map
between the pairs $Y=S\times A/(\iota,j)$ and $\check{Y}=
\check S\times A/(\check\iota,j)$ will take the form 
$$\check\phi:T_{Y}\rightarrow D_{\check{Y}}.$$
We start with 
$$B+i\omega=(B_1+i\omega_1,B_2+i\omega_2)\in T_{Y}\subseteq H^1(S^3,
R^1\bar\pi_*\C).$$ 
The map $\check{\bar \pi}:\check{Y}\rightarrow S^3$ is obtained
by dualizing $\bar\pi:Y\rightarrow S^3$, and hence by Poincar\'e
duality, $R^2\check{\bar\pi}_*\C\cong R^1\bar\pi_*\C$. Thus we obtain
an element $B+i\omega\in (F_2/F_1)\otimes_{\Q} \C$, the $E_2^{1,2}$
term of the Leray spectral sequence for $\check{\bar\pi}$. We claim that
there is a unique lift $\alpha$ of $B+i\omega$ to $F_2\otimes_{\Q} \C$
such that $E'\otimes s_x+\alpha\in D_{\check{Y}}$. This follows
immediately from the results for the K3 and elliptic cases. 
We then take $\check\phi(B+i\omega)$ to be $E'\otimes s_x+\alpha$. This
yields the desired mirror map.

In the above examples we have recovered previously known results, but
the construction used is of a form which (conjecturally) has far wider
application. In particular,
we use these examples to give a general conjecture as to how the mirror
map can be defined, which we hope will allow us to
construct mirrors and mirror maps whenever they exist.

As explained to the authors by Nigel Hitchin,
in order to specify a complex structure on a Calabi-Yau $n$-fold $X$,
we need to find a complex-valued $n$-form
$\Omega$ on $X$ which has the following properties:
\item{(1)}$\Omega$ is locally decomposable, i.e. locally $\Omega
=\theta_1\wedge\cdots\wedge \theta_n$ where $\theta_1,\ldots,\theta_n$
are 1-forms.
\item{(2)} $\Omega\wedge\bar\Omega$ is nowhere zero.
\item{(3)} $d\Omega=0$. (This condition tells us the almost complex
structure given by (1) and (2) is integrable.)

This data specifies a unique complex structure on $X$ in which $\Omega$
is a holomorphic $n$-form. Unfortunately, without a suitable Torelli
theorem, one does not know if the class $[\Omega]\in H^n(X,{\bf C})$
determines the complex structure uniquely.

We state a conjecture, which still needs further refinement to 
be totally precise. The results proved in this paper via the SYZ 
construction should be seen as evidence in favour of this conjecture.

\bigskip
{\it Conjecture.} Let $X$ be a Calabi-Yau threefold with compactified
moduli space $\overline{\M}_X$, with a large complex
structure limit point $P\in \overline{\M}_X$ (see [11] for a complete
definition of large complex structure limit point).
Then for some
open neighbourhood $U\subseteq\bar\M_X$ of $P$, and for general choice
of complex structure of $X$ in $U$ and general choice of K\"ahler class
$\omega$ on $X$, there is a special Lagrangian 3-torus fibration
$f:X\rightarrow B$ (with respect to the corresponding Ricci flat metric)
with fibre $T$ and with the following properties:
\item{(1)} $B$ is homeomorphic to $S^3$.
\item{(2)} $f$ has a special Lagrangian section $\sigma$.
\item{(3)} The Leray spectral sequence for $f:X\rightarrow B$
degenerates at the $E_2$ level and looks like
$$\matrix{
\Q \sigma&0&0&\Q\cr
0&E_2^{1,2}&H^4(X,\Q)&0\cr
0&H^2(X,\Q)&E_2^{2,1}&0\cr
\Q&0&0&\Q T\cr}$$
with $\dim E_2^{1,2}=\dim E_2^{2,1}=h^{1,2}(X)$, and the induced
filtration on $H^3(X,{\bf Q})$ coincides with the weight filtration
of the mixed Hodge structure associated with $P$.
\item{(4)} 
There is a natural compactification $f:\check X\rightarrow B$ of the
dual of $p:X_0\rightarrow B_0$ ($B_0\subseteq B$ the locus over
which $f$ is smooth, $X_0=f^{-1}(B_0)$), and the Leray spectral
sequence for $\check f:\check X\rightarrow B$ interchanges the first
and second rows of the spectral sequence for $f$.
\item{(5)} Let $B+i\omega$ be a general complexified K\"ahler class on $X$ for
``sufficiently large'' $\omega$. Then $B+i\omega\in H^2(X,{\bf C})$, which
is naturally isomorphic to the $E_2^{1,2}$ term of the spectral sequence
for $\check f:\check X\rightarrow B$, by (4). Let $[\check\sigma]$ be the
cohomology class of the canonical section of $\check f$. Then there
exists a lifting of $B+i\omega\in E_2^{1,2}$ to a class $\alpha
\in H^3(\check X,{\bf C})$, such that the cohomology class 
$[\check\sigma]+\alpha$ is represented by a locally decomposable
closed 3-form $\Omega$ on $\check X$ with $\Omega\wedge\bar\Omega$
nowhere zero. This induces a complex structure on $\check X$,
and this yields the mirror map, which is independent of the
original choices for complex structure and K\"ahler class.
\bigskip
{\it Remark 1.} We note that part (5) leaves open the question
of uniqueness of the lifting; there might be a number of choices,
and only the right choice yields the correct complex structure for
the mirror manifold. But in the examples considered here, namely
the K3 and Borcea-Voisin examples, the lifting is unique.
If the deformation family
of the complex structure on the mirror $\check X$
is already specified, then in fact,  from [4] the precise 
condition we need for uniqueness is that the full period
map is injective on some neighbourhood of large complex
structure limit; an elementary argument shows that it is
everywhere an immersion.
Another case
where this conjecture correctly puts the complex structure on the mirror
is for mirror symmetry of complex tori. 

{\it Remark 2.} The fact that the spectral sequence for $f:X\rightarrow B$
in the Borcea-Voisin examples does not look like the one given in the conjecture
should have to do with the fact we have used an $\omega$ which was on 
the boundary of the K\"ahler cone, as opposed to a general K\"ahler
class in the interior of the cone.

{\it Remark 3.} As a final example, we point out some evidence that
the SYZ construction is compatible with the Batyrev construction
for mirror symmetry (see [2]). This has also been observed by
D. Morrison and collaborators. If $\Delta$ is a reflexive polytope
of dimension $n$, giving rise to a toric variety $\P_{\Delta}$, then there
is a moment map $\mu:\P_{\Delta}\rightarrow\Delta$, whose general
fibre is a real $n$-torus. Now $\partial \Delta$ is homeomorphic to
$S^{n-1}$, and $X=\mu^{-1}(\partial\Delta)\subseteq\P_{\Delta}$ is
a union of toric divisors, and is a large complex structure limit
point in the family of Calabi-Yau $(n-1)$-folds in $\P_{\Delta}$.
The map $\mu:X\rightarrow \partial\Delta$ has general fibre a 
$T^{n-1}$. It should be possible, for small deformations
which smooth $X$, to also deform this torus fibration to yield a
special Lagrangian torus fibration on a smooth Calabi-Yau in the family.
\bigskip
\noindent\hd References\medskip

\rm\noindent 1.\ \quad Aspinwall, P.S., Morrison, D.R.:  String Theory on K3
surfaces.  
In:  Greene, B.R., Yau, S.-T. (eds.)  Essays on Mirror Manifolds II
(to appear).  Hong Kong, International Press 1996.\smallskip

\noindent 2.\ \quad Batyrev,  V.:  Dual Polyhedra and Mirror Symmetry for
Calabi-Yau Hypersurfaces in Toric Varieties.
J. Alg. Geom. {\bf 3}, 493--535 (1994).\smallskip

\noindent 3.\ \quad Borcea, C.:  K3 Surfaces with Involution and Mirror Pairs of
Calabi--Yau Manifolds.
In:  Greene, B.R., Yau, S.-T. (eds.)  Essays on Mirror Manifolds II
(to appear).  Hong Kong, International Press 1996.\smallskip

\noindent 4.\ \quad Bryant, R., Griffiths, P.: Some Observations on the
Infinitesimal Period Relations for Regular Threefolds with Trivial Canonical
Bundle.
In:  Artin, M., Tate, J. (eds.)  Arithmetic and Geometry, papers dedicated
to I.R. Shafarevich, Boston, Birkh\"auser, 1983, Vol. 2, 77--102.\smallskip

\noindent 5.\ \quad Cox, D., Zucker, S.:  Intersection Numbers of sections
of Elliptic Surfaces.
Inv. Math. {\bf 53}, 1-44 (1979).\smallskip

\noindent 6.\ \quad Dolgachev, I.V.:  Mirror Symmetry for Lattice Polarized K3
surfaces.  To appear, Duke eprint alg-geom/9502005.\smallskip

\noindent 7.\ \quad Dolgachev, I.V., Nikulin, V.V.:  Exceptional singularities of
V.I. Arnold and K3 surfaces.  Proc. USSR Topological Conf. in Minsk, 1977.
\smallskip

\noindent 8.\ \quad Harvey, R., Lawson, H.B. Jr.:  Calibrated Geometries.  Acta
Math. \bf 148\rm , 47-157 (1982).\smallskip

\noindent 9.\ \quad Hitchin, N.: Mirror Symmetry--- notes on Strominger et al.
Private Communication (1996).

\noindent 10.\ \quad McLean, R.C.:  Deformations of Calibrated Submanifolds. 
Texas A \& M University Preprint, 1996.\smallskip

\noindent 11.\quad Morrison, D.:  Compactifications of Moduli Spaces Inspired
by Mirror Symmetry. 
In: Journ\'ees de G\'eometrie Alg\'ebrique d'Orsay, Juillet 1992, Asterisque
\bf 218\rm , 243-271 (1993).

\noindent 12.\ \quad Morrison, D.R.: The Geometry Underlying Mirror Symmetry.
Preprint 1996.\smallskip

\noindent 13.\ \quad Nikulin, V.V.: Discrete reflection groups in Lobachevsky
spaces and algebraic surfaces.  Proc. ICM, Berkeley 1986, pp. 654-671.\smallskip

\noindent 14.\ \quad Pjatecki\u i-\u Sapiro, I.I., \u Safarevi\u c, I.R.:  
A Torelli Theorem
for algebraic surfaces of Type K3.  Math. USSR Izvestiya \bf 5\rm , 547-598
(1971).\smallskip

\noindent 15.\quad Pinkham, H.:  Singularit\'es exceptionelles, la dualit\'e
\'etrange d'Arnold et les surfaces K3.  C.R. Acad. Sci. Paris, Ser A, \bf
284\rm , 615-618 (1977).\smallskip

\noindent 16.\quad Strominger, A., Yau, S.-T., Zaslow, E.:  Mirror Symmetry is
T-Duality.  To appear, eprint hep-th/9606040.\smallskip

\noindent 17.\quad Voisin, C.:  Miroirs et involutions sur les surfaces K3. 
In: Journ\'ees de G\'eometrie Alg\'ebrique d'Orsay, Juillet 1992, Asterisque
\bf 218\rm , 273-322 (1993).

\hfill\vfill
\end